\documentclass{elsarticle}

\usepackage{graphicx,color,epsf,bm}
\usepackage{amsmath,amsfonts,amssymb,amscd,bm}
\usepackage[colorlinks=false, urlcolor=blue]{hyperref}
\usepackage[english]{babel}
\usepackage{graphics}
\usepackage{subcaption}
\usepackage{epstopdf}
\usepackage{afterpage}
\usepackage{natbib}

\begin{document}

\begin{frontmatter}

\title{Isogeometric Simulation of Lorentz Detuning 
in Superconducting Accelerator Cavities}

\address[POLI_MOX]{MOX Modeling and Scientific Computing, \\ 
                   Dipartimento di Matematica,
                   Politecnico di Milano, \\
                   Piazza L. da Vinci 32, 20133 Milano Italy}

\address[CEN]{CEN Centro Europeo di Nanomedicina,\\
              Piazza L. da Vinci 32, 20133 Milano Italy}

\address[GS]{Graduate School of Computational Engineering,\\ 
             Technische Universit\"at Darmstadt,\\
             Dolivostra\ss e 15 D-64293 Darmstadt, Germany}

\address[TEMF]{Institut f\"ur Theorie Elektromagnetischer Felder,\\
               Technische Universit\"at Darmstadt, \\
               Schlo\ss gartenstr. 8 64289 Darmstadt, Germany}

\author[GS,TEMF,POLI_MOX]{Jacopo Corno\corref{cor1}}
\ead{corno@gsc.tu-darmstadt.de}
\cortext[cor1]{Corresponding author}

\author[POLI_MOX,CEN]{Carlo de Falco}
\author[TEMF]{Herbert De Gersem}
\author[GS,TEMF]{Sebastian Sch\"{o}ps}

\begin{abstract}

Cavities in linear accelerators suffer from eigenfrequency shifts due to mechanical deformation caused by the electromagnetic radiation pressure, a phenomenon known as Lorentz detuning. Estimating the frequency shift up to the needed accuracy by means of standard Finite Element Methods, is a complex task due to the non exact representation of the geometry and due to the necessity for mesh refinement when using low order basis functions.
In this paper, we use Isogeometric Analysis for discretising both mechanical deformations and electromagnetic fields in a coupled multiphysics simulation approach. The combined high-order approximation of both leads to high accuracies at a substantially lower computational cost.

\end{abstract}

\begin{keyword}
Particle Accelerators; Superconducting Cavities; Isogeometric Analysis
\end{keyword}

\end{frontmatter}

\section{Introduction}

Controlling the resonant frequency of cavity eigenmodes in a particle accelerator is crucial in order to guarantee the synchronization of the electromagnetic wave and the particle bunches. Such frequency is determined essentially by the geometry of the cavity walls, which is therefore a critical parameter for the design of the cavity. The high-energy electromagnetic field inside the cavity exerts a radiation pressure on the walls, which causes a mechanical deformation of the geometry. Albeit small, this deformation may lead to a significant shift of the resonant frequency. This effect, known as \emph{Lorentz detuning}~\cite{Devanz_2002aa,Gassot_2002aa,Zaplatin_2006aa,Delayen_2003ab}, needs to be predicted with high precision in order to achieve a robust cavity design.

Standard Finite Element Methods (FEM) may require an extremely high level of mesh refinement to achieve sufficient accuracy when evaluating  Lorentz detuning, due to inaccuracies when approximating the deformed {and undeformed} cavity walls in the FEM mesh and due to the limited accuracy of typical low-order FEM basis functions. In this work, we propose a simulation strategy based on Isogeometric Analysis (IGA)~\cite{Cottrell2009} which allows an exact representation of the geometry and {the direct application of the computed deformation to the starting geometry, without any further approximation}. Finally it offers the possibility to accurately approximate the electromagnetic fields using high-order elements~\cite{Buffa2010}.

The outline of this paper is as follows: first we introduce the coupled electro\-magnetic-mechanical model describing Lorentz detuning. In the subsequent section Isogeometric Analysis is introduced along with an overview on the particular discretization used for Maxwell's equations. Finally we present the results obtained for the standard cylindrical test case and for the TESLA cavity geometry \cite{Tesla2000}.

\begin{figure}
	\begin{center}
	\includegraphics[width=.4\textwidth]{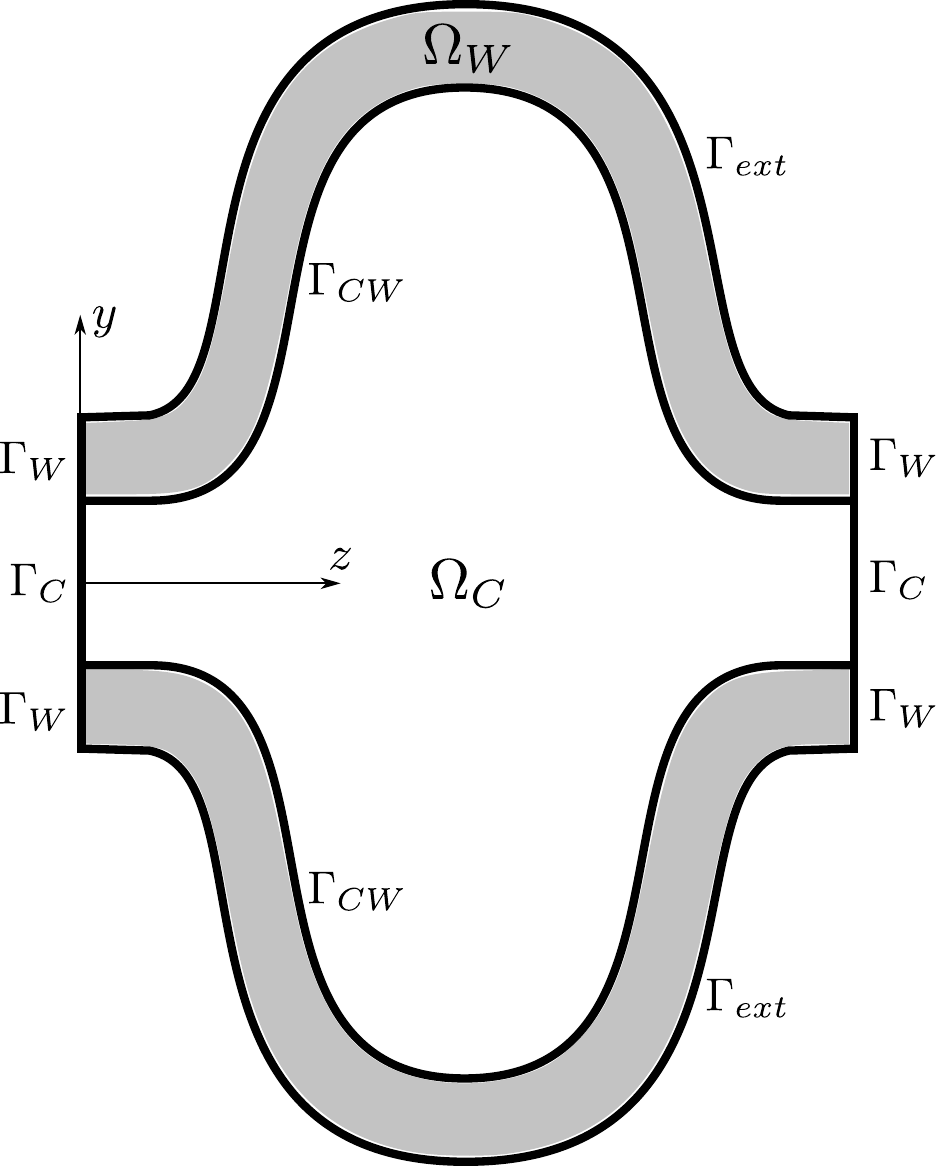}
	\par
	\caption{2D cut of the 3D computational domain for simulating Lorentz detuning in one cell of the TESLA cavity~\cite{Tesla2000} {(not to scale)} and labels for the domains and the boundaries ($yz$ section). The full cell is the result of a revolution around the $z$~axis.\label{fig:Domains}}
	\end{center}
\end{figure}

\section{Multi-physics Model for Lorentz Detuning}\label{sec:scheme}

Consider a one cell cavity geometry as the one depicted in Fig.~\ref{fig:Domains}. Let the two disjoint open domains with Lipschitz continuous boundaries $\Omega_{W} \subseteq \mathbb{R}^{3}$ and $\Omega_{C} \subseteq \mathbb{R}^{3}$ represent the cavity walls and the interior of the cavity, respectively. Let	 $\Gamma_{CW}=\overline{\Omega}_{C}\cap\overline{\Omega}_{W}$ denote the interface between the two domains. To evaluate the frequency shift, it is necessary to solve Maxwell's eigenproblem inside the undeformed and deformed cavity and an elasticity problem in the cavity walls.  We employ linear elasticity theory since the deformations are very small. The radiation pressure on the common interface $\Gamma_{CW}$ introduces a coupling between the two problems~\cite{Schreiber_2006aa}. The calculation steps are as follows:

\emph{Step 1.} Solve Maxwell's eigenproblem in $\Omega_{C}$:
\begin{subequations}\label{eq:Model_eigen}
    \begin{equation}\label{eq:Model_eigen_1}
	\nabla\times\left(\dfrac{1}{\mu_{0}}\nabla\times\mathbf{E}\right)=\omega_{0}^{2}\epsilon_{0}\mathbf{E} 
	\qquad\text{in }\Omega_{C}
	\end{equation}
	with the boundary conditions
    \begin{equation}\label{eq:Model_eigen_2}
    \begin{cases}
	\mathbf{E}\times\mathbf{n}_c=0 & \text{on }\Gamma_{CW}\\[2mm]
	\left(\dfrac{1}{\mu_{0}}\nabla\times\mathbf{E}\right)\times\mathbf{n}_c=0 & \text{on }\Gamma_{C}
	\end{cases}
    \end{equation}
\end{subequations}
where $\mu_0$ and $\epsilon_0$ are the permeability and permittivity of vacuum and $\mathbf{n}_c$ is the outward unit normal to $\Omega_C$. We assume time-harmonic fields with $\mathbf{E}$ a phasor {given in terms of peak values}. As cavity walls are often composed of a superconducting material, {\it e.g.} niobium, in order to reduce losses, they are assumed here to behave as a perfectly conducting boundary. At the two irises $\Gamma_{C}$, a Neumann condition is enforced, which is a common approximation corresponding to assuming the cell to be one of an infinite chain of cells. The eigenmode solution delivers a number of eigenfunction-eigenvalue  couplets, corresponding to the possible modes within the cavity. The accelerating mode of interest is the first transverse magnetic mode ($TM_{010}$). Let $\mathbf{E}_{0}$ be the computed electric field and $\omega_{0}^{2}$ the corresponding eigenvalue, then $f_{0}=\dfrac{\omega_{0}}{2\pi}$ is the resonant frequency for the accelerating eigenmode in the undeformed geometry.

\emph{Step 2.} Compute the magnetic field $\mathbf{H}_{0}$ for the first accelerating
eigenmode as
\begin{equation}
	\mathbf{H}_{0}=\dfrac{i}{\omega_{0}\mu_{0}}\nabla\times\mathbf{E}_0.
\end{equation}

The accelerating mode exerts on the cavity walls a radiation pressure with one component at $0$ frequency and one component at frequency $2\ f_{0}$. In practice, the latter can be neglected and the radiation pressure on $\Gamma_{CW}$ {is approximated by a time-constant value that} may be expressed as
\begin{eqnarray}
	\nonumber p&=& -\frac{1}{4}\epsilon_{0}\left(\mathbf{E}_0\:\mathbf{n}_c\right)\cdot\left(\mathbf{E}_0^{*}\;\mathbf{n}_c\right) \\
	&& +\frac{1}{4}\mu_{0}\left(\mathbf{H}_{0}\times\mathbf{n}_c\right)\cdot\left(\mathbf{H}^{*}_{0}\times\mathbf{n}_c\right)
\end{eqnarray}
{where $\mathbf{E}_0$ and $\mathbf{H}_0$ are peak values and $(\cdot)^{*}$ denotes the complex conjugate.}

\emph{Step 3.} Solve the following linear elasticity problem in the walls domain $\Omega_{W}$
\begin{subequations}\label{eq:Model_displ}
\begin{equation}\label{eq:Model_displ_1}
	\nabla\cdot\left(2\eta\nabla^{\left(S\right)}\mathbf{u}+\lambda\mathbf{I}\nabla\cdot\mathbf{u}\right)=0   
    \qquad \text{in }\Omega_{W}
	\end{equation}
	with boundary conditions
    \begin{equation}\label{eq:Model_displ_2}\small
	\begin{cases}
	\mathbf{u}=0 & \text{on }\Gamma_{W}\\[1mm]
	\left(2\eta\nabla^{\left(S\right)}\mathbf{u}+\lambda\mathbf{I}\nabla\cdot\mathbf{u}\right)\:\mathbf{n}_w=p\;\mathbf{n}_w & \text{on }\Gamma_{CW}\\[1mm]
	\left(2\eta\nabla^{\left(S\right)}\mathbf{u}+\lambda\mathbf{I}\nabla\cdot\mathbf{u}\right)\:\mathbf{n}_w=0 & \text{on }\Gamma_{ext}
	\end{cases}
	\end{equation}
\end{subequations}
for the displacement $\mathbf{u}$. In~\eqref{eq:Model_displ} we denote by $\nabla^{\left(S\right)}$ the symmetric gradient, while $\eta$ and $\lambda$ are the Lam\'{e} parameters of the wall constituent material and $\mathbf{n}_w$ is the outward unit normal to $\Omega_W$. On $\Gamma_{CW}$ the radiation pressure $p$ is applied.

\emph{Step 4.} Let the deformed walls domain $\Omega_{W}^{'}$ be defined as
\begin{equation}\label{eq:deformed_walls}
	\Omega_{W}^{'}\equiv\left\{\mathbf{x}+\mathbf{u}\left(\mathbf{x}\right),\,\mathbf{x}\in\Omega_{W}\right\},
\end{equation}
and the deformed cavity boundary $\Gamma_{CW}^{'}$ as
\begin{equation}
	\Gamma_{CW}^{'}\equiv\left\{\mathbf{x}+\mathbf{u}\left(\mathbf{x}\right),\,\mathbf{x}\in\Gamma_{CW}\right\}.
\end{equation}
Furthermore, let $\Omega_{C}^{'}$ denote the domain enclosed by $\Gamma_{CW}^{'}$
and $\Gamma_{C}$.

\emph{Step 5.} Solve Maxwell's eigenproblem in $\Omega_{C}^{'}$:

$$
	\nabla\times\left(\dfrac{1}{\mu_{0}}\nabla\times\mathbf{E}^{'}\right)=\left(\omega_{0}^{'}\right)^{2}\epsilon_{0}\mathbf{E}^{'} \qquad \text{in }\Omega_{C}^{'}
$$
with the boundary conditions
	\[
	\begin{cases}
	\mathbf{E}^{'}\times\mathbf{n}^{'}_c=0 & \text{on }\Gamma_{CW}^{'}\\[2mm]
	\left(\dfrac{1}{\mu_{0}}\nabla\times\mathbf{E}^{'}\right)\times\mathbf{n}^{'}_c=0 & \text{on }\Gamma_{C}^{'}
	\end{cases}
	\]
and let $\left(\left(\omega_{0}^{'}\right)^{2},\mathbf{E}_{0}^{'}\right)$ denote the accelerating eigenmode. The shifted frequency is finally obtained as
\[
	f_{0}^{'}=\dfrac{\omega_{0}^{'}}{2\pi}
\]
and the frequency shift due to Lorentz detuning as
\begin{equation}
	\Delta f_{0}=\left|f_{0}-f_{0}^{'}\right|.
\end{equation}

This procedure can be carried out iteratively if necessary.

\section{Numerical discretization}
Isogeometric Analysis (IGA) was born, less than a decade ago~\cite{Cottrell2005}, with the goal of \emph{bridging the gap between Computer Aided Design (CAD) and Finite Element Method (FEM)}. The main distinctive feature of IGA is that CAD geometries, commonly defined in terms of Non-Uniform Rational B-splines (NURBS), are represented exactly throughout the analysis, regardless of the level of mesh refinement, while in standard FEM the computational domain needs to be remeshed when performing \emph{h}-refinement and its geometry approaches the exact one only in the limit of vanishing mesh size \emph{h}.

Moreover, in addition to \emph{h}-refinement and \emph{p}-refinement, \emph{k}-refinement~\cite{Cottrell2009} was introduced as a combination of degree elevation and mesh refinement, yielding approximation spaces with higher regularity properties. \emph{k}-refinement has the advantage of not increasing the number of degrees of freedom of the problem, but produces matrices with larger bandwidth.

The particular IGA scheme adopted in this work takes advantage of the benefits of different approaches for each of the different physical subproblems being considered. The computational domains $\Omega_W$ and $\Omega_C$ are both defined via geometric mappings constructed in terms of NURBS basis functions. In solving the mechanical subproblem (\ref{eq:Model_displ}) an isoparametric approach is adopted so that the computed (discrete) displacement is defined in terms of the same NURBS basis and therefore the domain deformation (\ref{eq:deformed_walls}) is treated in a straight-forward way by a simple displacement of the control-points. In solving the Maxwell sub-problem (\ref{eq:Model_eigen}), on the other hand, the isoparametric approach is abandoned in favour of the choice of a solution space comprised of (push-forwards of) suitable B-spline functions which guarantees an $H\left(\mathrm{curl}\right)$ conforming, and therefore spectrally accurate, approximation of the field, as shown in~\cite{Buffa2010}. These concepts are explained in more detail below after introducing the required notation for NURBS and B-spline spaces.
\subsection{B-spline and NURBS functions}\label{sub:B-spline}
\begin{figure}[t]
  \centering
  \includegraphics[width=.9\textwidth]{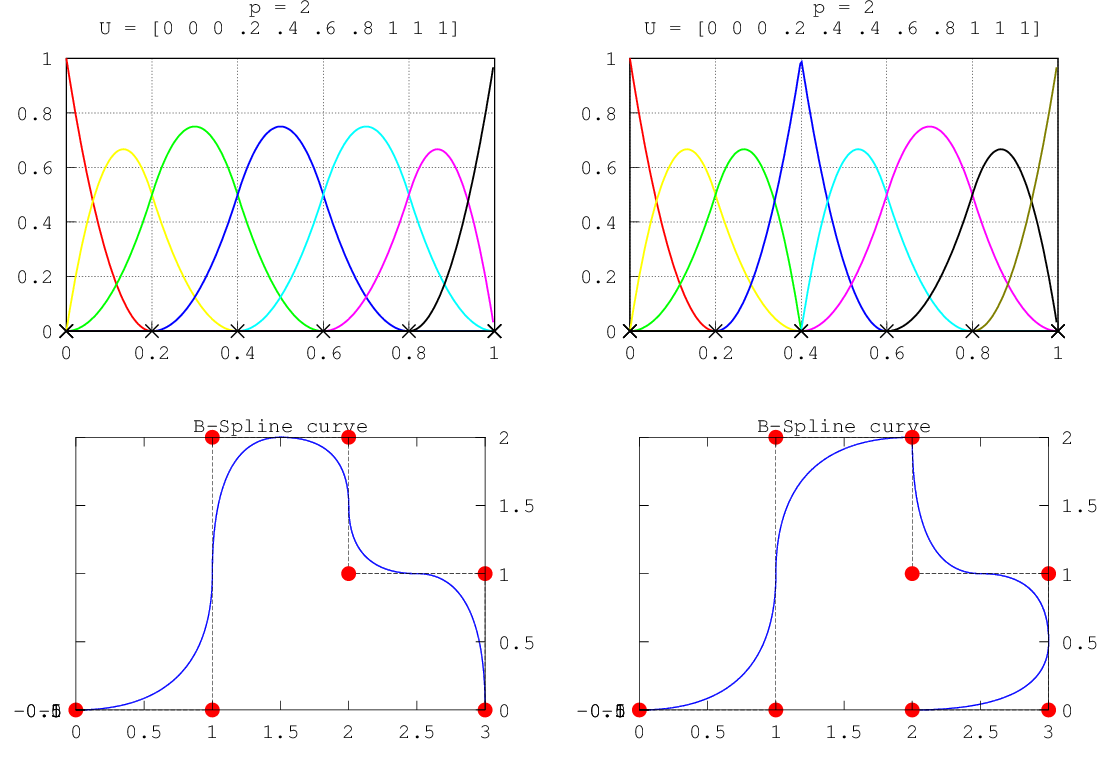}
  \caption{B-spline curves with different knot vectors: the multiplicity of the knot effects the regularity of the curve.}
\end{figure}
%
A B-spline geometrical entity is the result of the transformation through an appropriate mapping of a \emph{reference domain}. In one dimension, the reference domain is typically the interval $\left[0,1\right]$ which is then subdivided by a knot vector
\[
	\Xi=\left[\xi_0,\xi_1,\ldots,\xi_{n+p}\right]
\]
where $\xi_i\in\mathbb{R}\cap\left[0,1\right]$ is the $i$-th knot, $p$ is the polynomial degree ($p+1$ is the order) and $n$ is the number of basis functions used to build the B-spline curve. The knots divide the parameter space into elements. The element boundaries in the physical space are the images of the knots under the B-spline mapping.
Knot vectors can be uniform, if the knots are equally spaced, or non-uniform otherwise. Knots can be repeated and, by changing the multiplicity of a knot, we can change the level of continuity of the curve: basis functions of order $p$ have $p-r_i$ continuous derivatives across each knot $\xi_i$, where $r_i$ is the multiplicity of the $i$-th knot. In the particular case of a knot repeated exactly $r_i=p+1$ times, the basis is interpolatory at the knot $\xi_i$. A knot vector is said to be open if its first and last knots are repeated $p+1$ times (i.e. the curve is interpolatory at its ends). Below, we will always assume to be dealing with open knot vectors.

B-spline basis functions are defined by the Cox-De Boor recurrence formula:
\[
	B_{i,0}\left(\hat{x}\right)=
	\begin{cases}
		1 & \mathrm{if}\:\xi_i\leq\hat{x}\leq\xi_{i+1}\\
		0 & \mathrm{otherwise}
	\end{cases}
\]
\begin{align}\label{eq:Cox-deBoor}
  \begin{split}
	  B_{i,p}\left(\hat{x}\right) = & \frac{\hat{x}-\xi_{i}}{\xi_{i+p}-\xi_{i}}B_{i,p-1}\left(\hat{x}\right)+\\
	  & \frac{\xi_{i+p+1}-\hat{x}}{\xi_{i+p+1}-\xi_{i+1}}B_{i+1,p-1}\left(\hat{x}\right)
  \end{split}
\end{align}
with $i=0, \ldots, n-1$. We will denote the space spanned by the $n$ functions $B_{0,p},\ldots,B_{n-1,p}$ by $S_{\boldsymbol\alpha}^{p}\left(\Xi\right)$, with $\boldsymbol\alpha=\left\{ \alpha_0,\ldots,\alpha_{n+p+1}\right\}$ and $\alpha_{i}=p-r_{i}$, where $r_{i}$ is the multiplicity of the $\mathit{i}$-th knot.

B-spline curves are built taking a linear combination of B-spline basis functions and defining a set of control points. In particular, given $n$ basis functions $B_{i,p}$ and $n$ control points $\mathbf{P}_{i}\in\mathbb{R}^{d}$, $i=0,1,\ldots,n-1$,
a piecewise polynomial B-spline curve is defined by the following:
\begin{equation}\label{eq:BsplineCurve}
	\mathbf{C}\left(\hat{x}\right)=\sum_{i=0}^{n-1}B_{i,p}\left(\hat{x}\right)\mathbf{P}_{i}.
\end{equation}

The concepts presented until now can be easily extended to B-spline surfaces and volumes using a tensor product approach. For instance in the 3D case, given the knot vectors $\Xi_d$, the degrees $p_d$ and the number of basis funtions $n_d$ (with $d=1,2,3$), the B-spline trivariate basis functions are defined as
\begin{equation}
	B_{\mathbf{i}}^\mathbf{p}\left(\hat{\mathbf{x}}\right)=B_{i_1,p_1}\left(\hat{x}\right) B_{i_2,p_2}\left(\hat{y}\right) B_{i_3,p_3}\left(\hat{z}\right),
\end{equation}
where $\mathbf{p}=\left(p_1,p_2,p_3\right)$ and $\mathbf{i}=\left(i_1,i_2,i_3\right)$ is a multi-index in the set
\[
	\mathcal{I} = \left\{\mathbf{i}=\left(i_1,i_2,i_3\right) : 0\leq i_d \leq n_d-1\right\}.
\]
Given the regularities $\boldsymbol\alpha_1$, $\boldsymbol\alpha_2$, and $\boldsymbol\alpha_3$, we will refer to this space of B-splines as $S_{\boldsymbol\alpha_{1},\boldsymbol\alpha_{2},\boldsymbol\alpha_{3}}^{\mathbf{p}}$.

Starting from the Cox-de Boor formula given in (\ref{eq:Cox-deBoor}), we can define the rational basis functions $N_{\mathbf{i}}^{\mathbf{p}}\left(\hat{\mathbf{x}}\right)$:
\begin{equation}
	N_{\mathbf{i}}^{\mathbf{p}}\left(\hat{\mathbf{x}}\right)=\dfrac{B_{\mathbf{i}}^{\mathbf{p}}\left(\hat{\mathbf{x}}\right)w_{\mathbf{i}}}{\sum_{\mathbf{j}\in\mathcal{I}}B_{\mathbf{j}}^{\mathbf{p}}\left(\hat{\mathbf{x}}\right)w_{\mathbf{j}}}
\end{equation}
where we assume $w_{\mathbf{i}}>0$ for all $\mathbf{i}$. We will denote the space of NURBS with $\mathcal{N}^{\mathbf{p}}$. A NURBS object is built in an analogous way to (\ref{eq:BsplineCurve}):
\begin{equation}\label{eq:NURBSCurve}
  \mathbf{C}\left(\hat{\mathbf{x}}\right)=\sum_{\mathbf{i}\in\mathcal{I}}N_{\mathbf{i}}^{\mathbf{p}}\left(\hat{\mathbf{x}}\right)\mathbf{P}_{\mathbf{i}}.
\end{equation}

With respect to B-spline, using NURBS, one can utilize both the control points and the weights to control the local shape: as $w_{\mathbf{i}}$ increases, the curve is pulled closer to the control point $\mathbf{P}_{\mathbf{i}}$, and viceversa. This allow the exact representation of important geometries, often used in CAD, such as conic sections.

\subsection{Linear elasticity problem}
The weak formulation of (\ref{eq:Model_displ}) is:

\textit{Find the displacement $\mathbf{u}\in \left(H_{0}^{1}\left(\Omega_{W}\right)\right)^{3}$ such that}
\begin{multline}
  \int_{\Omega_W}\left(2\eta\varepsilon\left(\mathbf{u}\right):\varepsilon\left(\mathbf{v}\right)+\lambda\left(\nabla\cdot\mathbf{u}\right)\left(\nabla\cdot\mathbf{v}\right)\right)\, d\mathbf{x}=\\
  \int_{\Gamma_{N}}p\mathbf{n}_w\cdot\mathbf{v}\, d\Gamma,\quad\forall\mathbf{v}\in \left(H_{0}^{1}\left(\Omega_{W}\right)\right)^{3}
\end{multline}
where $\varepsilon=1/2\left(\nabla\mathbf{u}+\nabla\mathbf{u}^T\right)$ is the small deformation strain tensor. In structural mechanics, it is very useful to invoke the isoparametric concept, such that the undeformed and deformed geometry belong to the same function spaces. This means that when the problem is solved using a higher order solution space, the order of the geometrical representation has to be elevated accordingly. In the \emph{k}-refinement approach this is achieved via appropriate knot-insertion so that the shape of the domain is not changed.

Let the walls domain $\Omega_W$, bounded and Lipschitz, be represented by a NURBS volume 
\begin{equation}\label{eq:transformationF_W}
	\Omega_W=\mathbf{F}_W\left(\hat{\Omega}_W\right) = \sum_{\mathbf{i}\in\mathcal{I}}N_{\mathbf{i}}^{\mathbf{p}}\left(\hat{\mathbf{x}}\right)\mathbf{P}_{\mathbf{i}}
\end{equation}
where $\mathbf{F}_W$ is a 3D mapping of the type introduced in (\ref{eq:NURBSCurve}) (smooth with an almost everywhere piecewise smooth inverse). The discrete space is the space $V_h$ obtained by the transformation through $\mathbf{F}_W$ of the same space $\mathcal{N}^{\mathbf{p}}$ that defines the geometry:
\begin{equation}
	V_h = \left\{v_h\in H_{0,\Gamma_D}^1 : v_h=\hat{v}_h\circ\mathbf{F}_W^{-1},\hat{v}_h\in \mathcal{N}^{\mathbf{p}}\right\}
\end{equation}

With this choice, the deformed geometry is elegantly obtained by adding the solution vector $\mathbf{u}$ to the control net of the initial NURBS domain
\begin{equation}
	\Omega_W^{'} = \mathbf{F}_W^{'}\left(\hat{\Omega}_W\right) = \sum_{\mathbf{i}\in\mathcal{I}}N_{\mathbf{i}}^{\mathbf{p}}\left(\hat{\mathbf{x}}\right)\left(\mathbf{P}_{\mathbf{i}}+\mathbf{u}_{\mathbf{i}}\right).
\end{equation}

\subsection{Electromagnetic cavity eigenproblem}
Let $\Omega_{C}\in\mathbb{R}^{3}$ be our bounded NURBS cavity domain. Using Green's integration by parts formula and the notion of $H_{0,\Gamma_{D}}\left(\mathrm{curl};\Omega_{C}\right)$ of functions with curl well defined in $L^2$ and vanishing trace on the boundary, a standard variational formulation of problem (\ref{eq:Model_eigen}) reads as follows~\cite{Buffa2010}:

\textit{Find $\omega\in\mathbb{R}$, and $\mathbf{E}\in H_{0,\Gamma_{D}}\left(\mathrm{curl};\Omega_{C}\right)$, with $\mathbf{E}\neq0$, such that}
\begin{multline}\label{eq:eigen_var_form}
	\int_{\Omega_C}\mu_{0}^{-1}\nabla\times\mathbf{E}\cdot\nabla\times\mathbf{w}\, d\mathbf{x} = \omega^{2}\int_{\Omega_C}\epsilon_{0}\mathbf{E}\cdot\mathbf{w}\, d\mathbf{x}\\
	\forall\mathbf{w}\in H_{0,\Gamma_{D}}\left(\mathrm{curl};\Omega_{C}\right).
\end{multline}
It is known that $\omega =0$ is the essential spectrum, and that its associated eigenspace has infinite dimension. All other eigenvalues form a diverging sequence with associated eigenspaces belonging to $H_{0,\Gamma_{D}}\left(\mathrm{curl};\Omega_{C}\right)\cap H\left(\mathrm{div_{0}};\Omega_{C}\right)$, where we denote with $H\left(\mathrm{div_{0}};\Omega_{C}\right)$ the space of function in $H\left(\mathrm{div};\Omega_{C}\right)$ with divergence equal to zero.

The functional spaces used for the variational formulation (\ref{eq:eigen_var_form})
have some special relations that are summarized through the well known de Rham diagram~\cite{Buffa2010}. In order to achieve a consistent approximation of Maxwell's eigenvalue problem, the discrete spaces have to satisfy an analogous relation.

Following~\cite{Buffa2010}, we define on the reference domain a vectorial B-spline space with differing degree for each component:
\begin{equation}
	S^{1} = S_{\boldsymbol\alpha_{1}-\mathbf{1},\boldsymbol\alpha_{2},\boldsymbol\alpha_{3}}^{p_{1}-1,p_{2},p_{3}}\times S_{\boldsymbol\alpha_{1},\boldsymbol\alpha_{2}-\mathbf{1},\boldsymbol\alpha_{3}}^{p_{1},p_{2}-1,p_{3}}\times S_{\boldsymbol\alpha_{1},\boldsymbol\alpha_{2},\boldsymbol\alpha_{3}-\mathbf{1}}^{p_{1},p_{2},p_{3}-1}
\end{equation}
{where $\boldsymbol\alpha_{i}-\mathbf{1}$ states that the regularity at each knot is decreased by one (since the corresponding degree is decreased).}

The final step is to define the finite dimensional spaces in the physical domain $\Omega_{C}$. Let $\mathbf{F}_C$ be the parametrization for our domain computed with the same hypothesis as given for (\ref{eq:transformationF_W}), then the discrete space on $\Omega_C$ is defined through a curl conforming mapping \cite{Monk2003}: 
\begin{equation}\label{eq:discrete_space_curl_preserv}
	X^1=\left\{ \left(D\mathbf{F}_C\right)^{-T}\left(\mathbf{w}\circ\mathbf{F}_C^{-1}\right),\mathbf{w}\in S^{1}\right\}
\end{equation}
where $D\mathbf{F}_C$ is, the Jacobian matrix of the parametrization. It has been proven \cite{Buffa2010} that this space has the approximation properties needed for the discretization of $H\left(\mathrm{curl}\right)$.
\subsection{Multipatch formulation}
\begin{figure}[b]
\begin{center}
  \includegraphics[width=.8\textwidth]{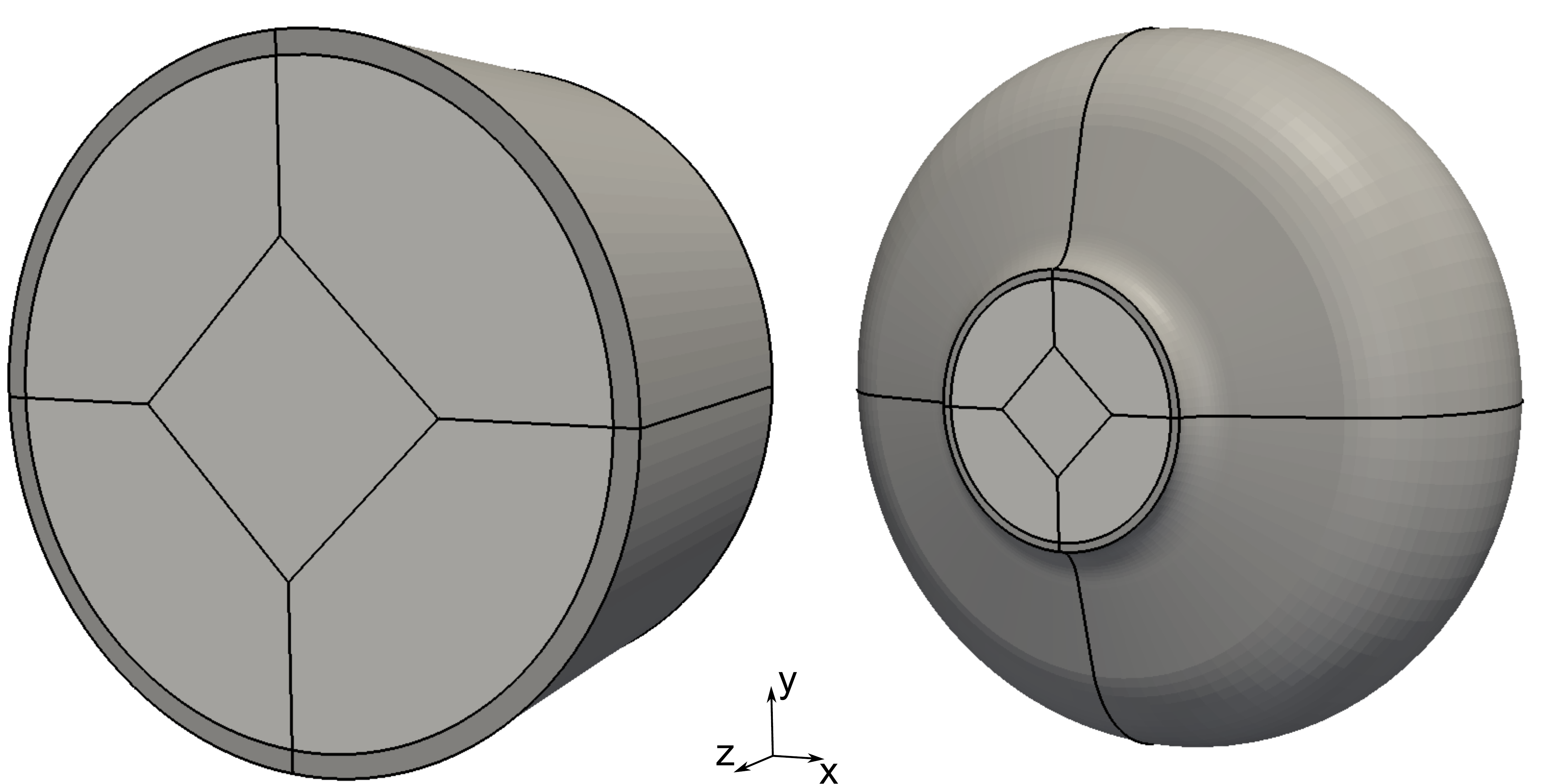}
  \caption{Patch subdivision for the pill-box cavity (left) and the TESLA cavity (right).}\label{fig:patches}
\end{center}
\end{figure}
%
In some situations, using a \emph{single patch} domain geometry definition as in (\ref{eq:transformationF_W}) is impossible or at least inconvenient.
For example, in parametrizing the geometries for both the cylindrical pill-box cavity 
and for the TESLA cavity, that are the focus of the present work, we have chosen to use 
a multipatch approach in order to avoid singularities in the geometrical mapping \cite{Veiga2012}. 
In other words the walls domain geometry for our problems is partitioned into $N_{w}$ subregions as
\begin{equation}
\begin{array}{l}
\overline{\Omega}_{W} \equiv \cup_{i=1}^{N_{w}}\; \overline{\Omega}_{W,i}\\[2mm]
\Omega_{W,i} \cap \Omega_{W,j} = \emptyset \qquad \forall i\neq j
\end{array}
\end{equation}
where each of the \emph{patches} consists of a smooth mapping with smooth inverse of the reference domain $\hat{\Omega}$
$$
\Omega_{W,i} \equiv \mathbf{F}_{W,i} (\hat{\Omega}),
$$
each of the mappings $\Omega_{W,i}$ being defined in terms of NURBS basis functions as in (\ref{eq:transformationF_W}).
We require that two neighbouring patches share one full face and we denote the interface by
$$
\Gamma_{W,ij} \equiv \overline{\Omega}_{W,i} \cap \overline{\Omega}_{W,j}.
$$
The resulting overall geometrical mapping is globally continuous but only piecewise smooth.
A similar partitioning and similar notation is used for the multipatch parametrization of the cavity domain, {\it i.e.}
\begin{equation}
\begin{array}{l}
\overline{\Omega}_{C} \equiv \cup_{i=1}^{N_{c}}\; \overline{\Omega}_{C,i}\\[2mm]
\Omega_{C,i} \cap \Omega_{C,j} = \emptyset \qquad \forall i\neq j,
\end{array}
\end{equation}
with
$$
\Omega_{C,i} \equiv \mathbf{F}_{C,i} (\hat{\Omega}),
$$
and
$$
\Gamma_{C,ij} \equiv \overline{\Omega}_{C,i} \cap \overline{\Omega}_{C,j}.
$$

In Fig.~\ref{fig:patches} the subdivisions for the two geometries being considered in this paper are depicted.
To extend the linear elasticity (\ref{eq:Model_displ}) and Maxwell (\ref{eq:Model_eigen}) problem to the new geometric setting, a substructuring approach is used.
For the problem (\ref{eq:Model_displ}) a new set of unknowns $\mathbf{u}_{i}$ is introduced, such that $\left.\mathbf{u}\right|_{\Omega_{W,i}} = \mathbf{u}_{i}$ and a problem similar to~(\ref{eq:Model_displ_1}) is set in each patch
$$
\nabla\cdot\left(2\eta\nabla^{\left(S\right)}\mathbf{u}_i+\lambda\mathbf{I}\nabla\cdot\mathbf{u}_i\right)=0   
    \qquad \text{in }\Omega_{W,i}
$$
and the overall problem~(\ref{eq:Model_displ}) is recovered by imposing continuity of the displacements and normal stresses at the patch interfaces
$$
	\mathbf{u}_i = \mathbf{u}_j \qquad \text{on } \Gamma_{W,ij}
$$
\vspace{-.7cm}
\begin{multline*}	
	\left(2\eta\nabla^{\left(S\right)}\mathbf{u}_i+\lambda\mathbf{I}\nabla\cdot\mathbf{u}_i\right)\cdot\mathbf{n}_{w,i} + \\
	\left(2\eta\nabla^{\left(S\right)}\mathbf{u}_j+\lambda\mathbf{I}\nabla\cdot\mathbf{u}_j\right)\cdot\mathbf{n}_{w,j} = 0 \; \text{on } \Gamma_{W,ij}.
\end{multline*}

Similarly for the Maxwell eigenproblem the unknowns $\mathbf{E}_i$ are introduced and the problem to be solved in each patch becomes
$$
\nabla\times\left(\dfrac{1}{\mu_{0}}\nabla\times\mathbf{E}_i\right)=\omega_{0}^{2}\epsilon_{0}\mathbf{E}_i \qquad \text{in }\Omega_{C,i}
$$
with the interface conditions
$$
	\mathbf{E}_i\times\mathbf{n}_{c,i} = \mathbf{E}_j\times\mathbf{n}_{c,j} \qquad \text{on } \Gamma_{C,ij}
$$
$$
	\left(\dfrac{1}{\mu_{0}}\nabla\times\mathbf{E}_i\right)\times\mathbf{n}_{c,i} +
 	\left(\dfrac{1}{\mu_{0}}\nabla\times\mathbf{E}_j\right)\times\mathbf{n}_{c,j} = 0 \;\text{on }\Gamma_{C,ij}.
$$

\begin{figure}
  \centerline{\includegraphics[width=.65\columnwidth]{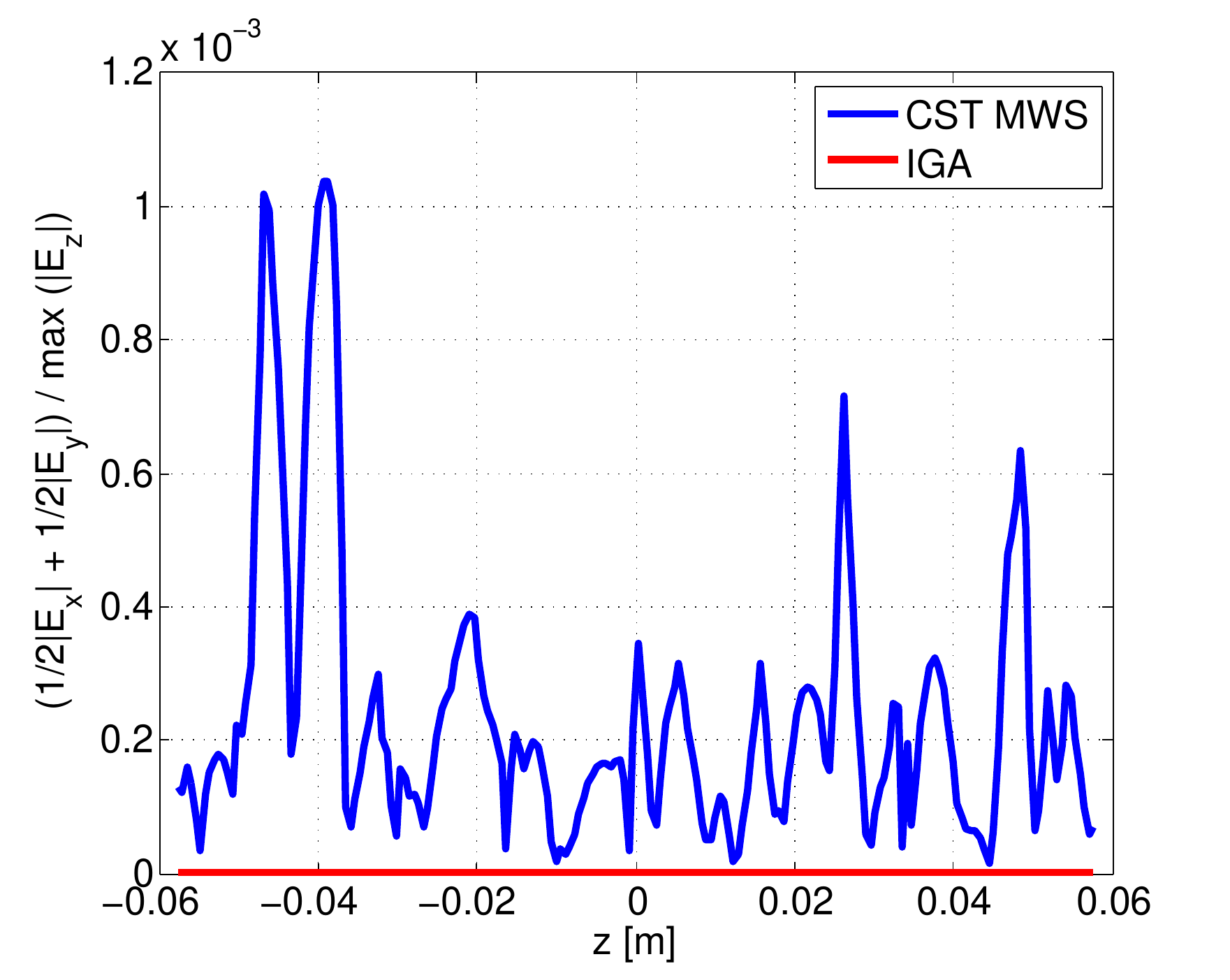}}
  \caption{Oscillations in the transverse component of the electric field along the axis of the 1-cell TESLA cavity. The correct accelerating field should have $E_x = E_y = 0$ and only longitudinal component. The FEM method suffers from oscillations due to non axis-aligned elements, while the IGA solution is precisely determined. Both results are obtained with a second order approximation.}\label{fig:oscillations}
\end{figure}

With respect to standard FEM, where only the tangential component of the computed solution is continuous across the elements boundaries, given IGA high regularity properties it is possible to achieve solutions with higher smoothness (up to $C^{p-1}$) within each patch. Only across the patch interfaces the regularity is reduced to $C^0$. The patches have been created in such a way that there are no interfaces across the length ($z$ direction) of the cavity and this is of great interest since it leads to smooth solutions particularly along the $z$-axis of the cavity, where the particle bunches travel and thus high precision is required. Classical FEM cavity simulations on tetrahedra may not achieve sufficient precision since the solution is often affected by undesired oscillations due to the discontinuities across the elements and not axis-aligned edges (see Fig.~\ref{fig:oscillations}). Usually this problem is solved by using a huge number of tetrahedra, symmetric or hybrid meshes, e.g., with hexahedra along the axis~\cite{Gjonaj2009}. By using an Isogeometric mesh it is possible to completely avoid the problem in an easy and computationally inexpensive way.

\section{Results}

\begin{figure}
\begin{subfigure}{.5\textwidth}
	\centerline{\includegraphics[width=\textwidth]{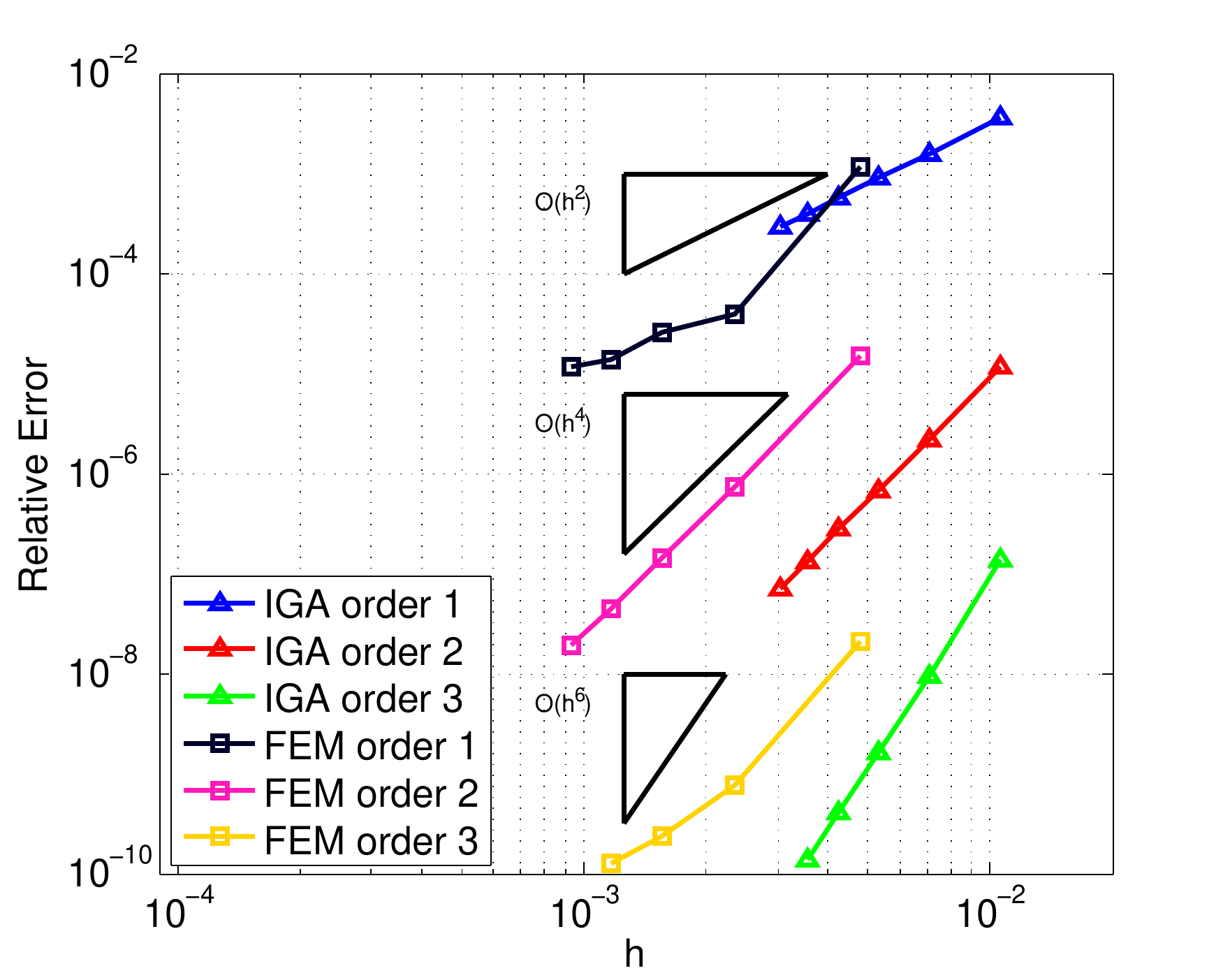}}
  \caption{Comparison of IGA proposed method and CST implementation with sensitivity analysis}\label{fig:Convergence}
\end{subfigure}
~
\begin{subfigure}{.5\textwidth}
  \centerline{\includegraphics[width=\columnwidth]{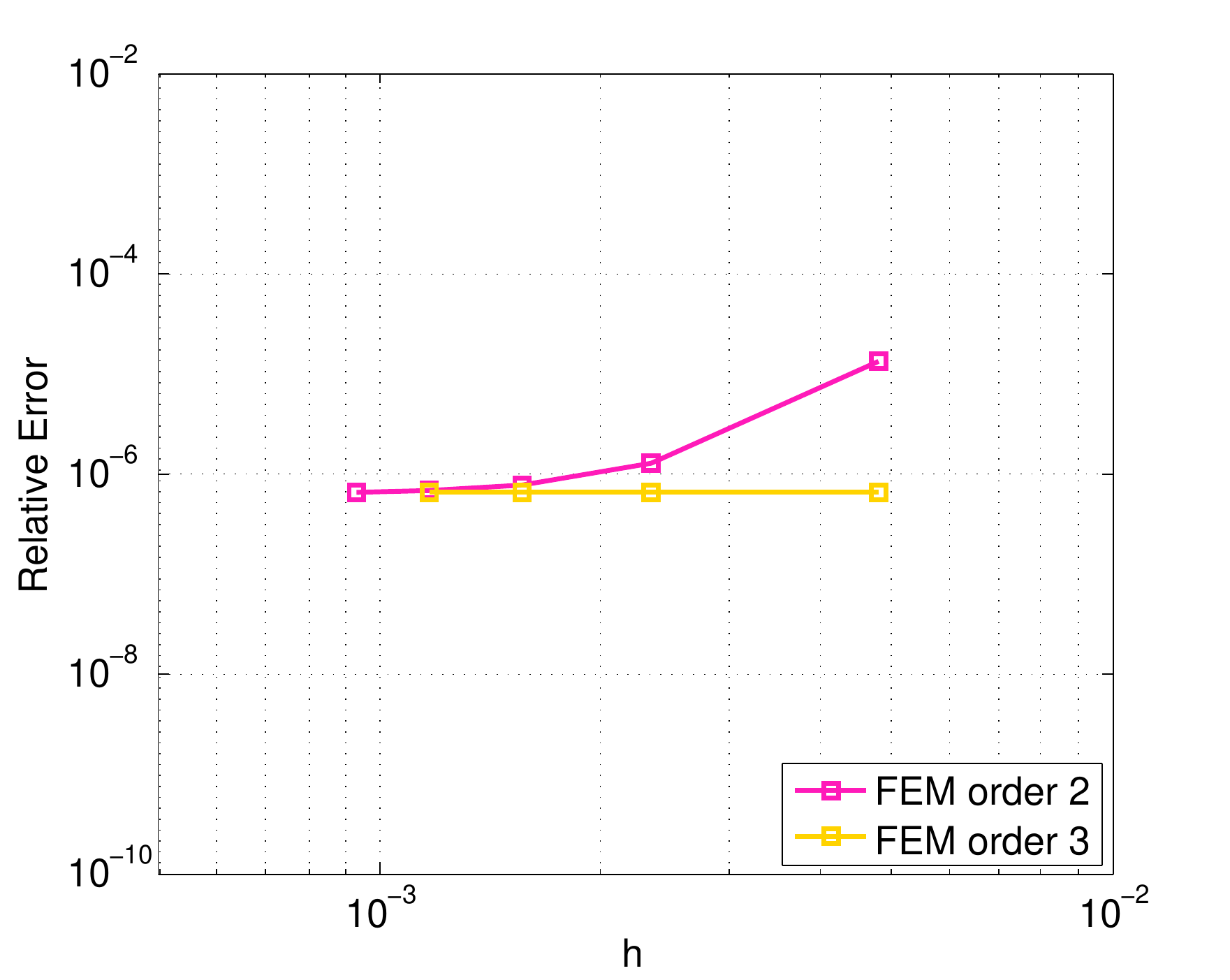}}
  \caption{Convergence of CST results for the resonant frequency with the algorithm proposed in Section~\ref{sec:scheme}.}\label{fig:CST_aba}
  \end{subfigure}
\caption{Convergence of the eigenfrequency for the deformed pill-box cavity (design parameters: $R=35$ mm, $L=100$ mm, exact frequency $f_0^{'}=3.278292919$ GHz). The IGA simulation is performed following the steps described in Section~\ref{sec:scheme}. The FEM results are obtained using the commercial software CST.}
\end{figure}

\begin{figure}
	\begin{center}
		\includegraphics[height=0.3\textheight]{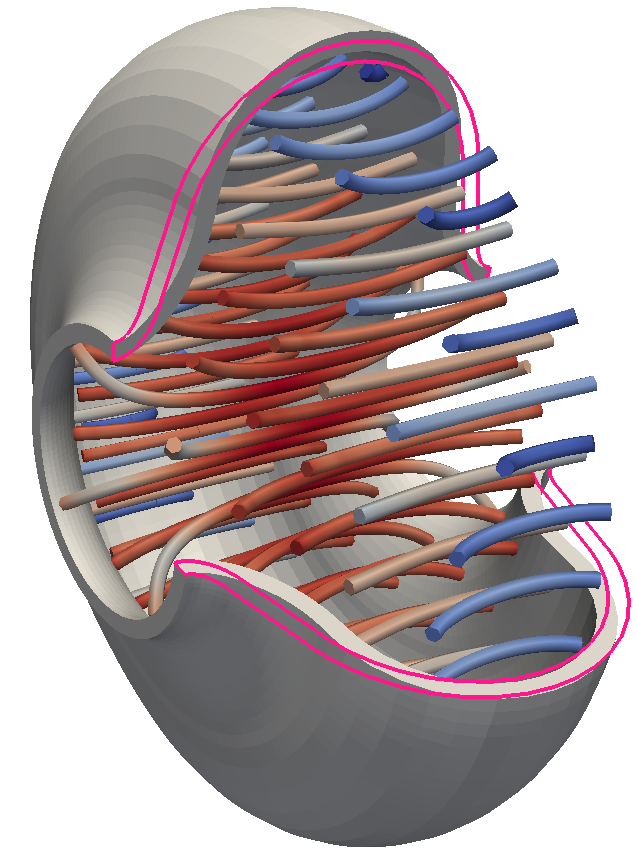} \;\hspace{1cm}\;
		\includegraphics[angle=90,height=0.3\textheight]{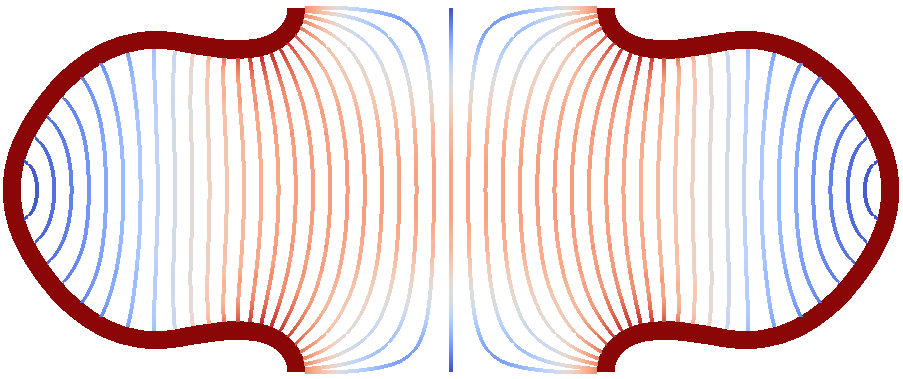}
	\caption{On the left: undeformed geometry (pink) and deformed geometry (grey, amplified by	a factor $5\cdot10^5$) for the 1-cell TESLA cavity. On the right: deformed accelerating cavity mode.\label{fig:Displ3d}}
	\end{center}
\end{figure}

The implementation of the discretization scheme just introduced has been done in MATLAB~\cite{MATLAB} and Octave~\cite{Octave} using GeoPDEs~\cite{GeoPDEs}.
Its applicability for cavity simulation has been verified by using a pill-box cavity with known closed form solution~\cite{Jackson1962} (excluding the bases of the cylinder in the mechanical simulation, i.e. the resulting deformation is only radial). The steps illustrated in Section~\ref{sec:scheme} have been applied to the first transverse magnetic ($TM_{010}$) mode in the cavity and the corresponding detuning has been computed. {The eigenvalue problem has been solved using the ARPACK library and the Implicitly Restarted Lanczos Method (IRLM), a variant of the Arnoldi/Lanczos process with the Implicitly Shifted QR technique that is suitable for large sparse matrices \cite{Lehoucq_1997aa}.} The new value of the frequency has been compared with the exact solution given by the theory while increasing the mesh resolution for a given polynomial order (see Fig.~\ref{fig:Convergence}). Of particular relevance is the fact that the multiphysical coupling does not decrease the optimal  convergence rates for the eigenvalue problem.

{As a comparison, a similar procedure was performed in the proprietary electromagnetic field simulation software CST STUDIO SUITE\textsuperscript{\textregistered}~\cite{CST} that is the quasi-standard for cavity simulation. The eigenproblem in the cavity is solved in CST MICROWAVE STUDIO\textsuperscript{\textregistered} (MWS) using the FE eigenvalue solver and the Lorentz forces are exported to CST MPHYSICS STUDIO\textsuperscript{\textregistered} (MPS) to compute the wall deformation. The information on the displacement is then imported once again in MWS, where the detuned frequency $f_0^{'}$ is estimated through a sensitivity analysis approach. The results are depicted in Fig.~\ref{fig:Convergence} along the IGA ones. The approach used in CST leads to a linearization of the problem but the method performs well since the deformations are very small. In addition to the simulation using sensitivity analysis, the proposed algorithm for IGA was implemented in CST. The results show that the level of accuracy reachable in this case for the resonating frequency is limited to $10^{-6}$ (see Figure~\ref{fig:CST_aba}).}

\begin{table}[t]
\begin{center}
\resizebox{\textwidth}{!}{%
\begin{tabular}{|c||r|r|r|r||r|r|r|r|}
			\hline 
& \multicolumn{4}{|c||}{\textbf{2$^\text{nd}$ order}} & \multicolumn{4}{c|}{\textbf{3$^\text{rd}$ order}}
 	\tabularnewline
& \multicolumn{2}{|c}{\textbf{IGA}} & \multicolumn{2}{c||}{\textbf{FEM}} &
 	\multicolumn{2}{|c}{\textbf{IGA}} & \multicolumn{2}{c|}{\textbf{FEM}}
 	\tabularnewline\hline
\textbf{Rel. error} & 
\multicolumn{1}{|c|}{$N_{\text{dof}}$} & \multicolumn{1}{|c|}{$t$ [s]} & \multicolumn{1}{|c|}{$N_{\text{dof}}$} & \multicolumn{1}{|c||}{$t$ [s]} & \multicolumn{1}{|c|}{$N_{\text{dof}}$} & \multicolumn{1}{|c|}{$t$ [s]} &
\multicolumn{1}{|c|}{$N_{\text{dof}}$} & \multicolumn{1}{|c|}{$t$ [s]}
\tabularnewline\hline \hline
1e-05 &  1540	&  0.2 &   5346 &   1.7 &       &      &      &\tabularnewline \hline
1e-06 &  9828 &  6.8 &  46266 &  21.1 &       &      &      &\tabularnewline \hline
1e-07 & 18304 & 14.8 & 158050 & 187.6 &       &      &      &\tabularnewline \hline
1e-08 & 47520 & 95.1 & 381036 & 843.4 &  4480 &   2.5 &  15618 &    5.8\tabularnewline \hline
1e-10 &       &      &        &       & 30628 &  91.7 & 135246 &  141.5\tabularnewline \hline
1e-11 &       &      &        &       & 97888 & 542.8 & 461937 & 1176.3\tabularnewline \hline
\end{tabular}
}
\caption{Number of DoFs required to compute the first accelerating mode in the pill-box cavity within a prescribed accuracy ($R=35$ mm and $L=100$ mm, $f_0=3.2783579381$ GHz). The IGA implementation was performed in GeoPDEs~\cite{GeoPDEs} while for the FEM simuation CST STUDIO SUITE~\cite{CST} was used (empty cells are due to unavailable FEM matrices). The times listed refer to the solution of the eigenvalue problem with ARPACK.\label{tab:PillboxFEMvsIGA}}
\end{center}
\end{table}

{In order to be able to fairly compare the two codes in terms of efficiency, a set of matrices, with increasing mesh resolution, was generated in CST for 2$^{\text{nd}}$ and 3$^{\text{rd}}$ order basis functions and exported to MATLAB. The same Arnoldi solver used for the IGA matrices was applied to solve the generalized eigenvalue problem for the FEM ones. In Table~\ref{tab:PillboxFEMvsIGA} we report the number of degrees of freedom required by the IGA and FEM methods to achieve a given level of accuracy, alongside with the time needed to solve the corresponding eigenvalue problem. Since the B-spline basis functions have a wider support, the IGA matrices are denser than their FEM counterparts. For example, given an Isogeometric matrix of dimension $50000$ approximately, the ratio of non zero elements over total number of elements is $2.4\cdot10^{-2}$, while for an analogous FEM matrix the ratio is $8.8\cdot10^{-4}$. However the accuracy-per-degree-of-freedom is higher when using Isogeometric Analysis and this leads to speed-ups up to 9 times (2$^{\text{nd}}$ order, error 1e-8) as shown in Table~\ref{tab:PillboxFEMvsIGA}.}

A second more realistic example is the 1-cell TESLA cavity \cite{Tesla2000} (see Fig.~\ref{fig:Domains}). The accelerating eigenmode of the TESLA cavity is the $TM_{010}$ mode at $1.3$~GHz. The frequencies for undeformed and deformed geometry are computed on six meshes with an increasing number of subdivisions (Table~\ref{tab:FrequencyShift}). In the last column of  Table~\ref{tab:FrequencyShift} we report the difference between the values of the frequency shift computed at two subsequent levels of refinement, which shows that six subdivisions, corresponding to about $110000$ DOFs, are sufficient to achieve an accuracy of about 1~Hz. In this last case, the total computational time (geometry creation, matrix construction and eigenvalue solver) is approximately 10-15 minutes.  In Fig.~\ref{fig:Displ3d}, the undeformed and deformed geometry are compared. The computed displacement is in the order of 1~nm~$\sim$~10~nm, which is in good accordance to results reported in literature \cite{Tesla2000}.

\begin{table}
\begin{center}
		\begin{tabular}{|c|c|c|c|c|c|}
			\hline 
			subs & $N_{el}$ & $N_{\text{dof}}$ & $f_0$ [GHz]& Shift [Hz]& variation  [Hz] \tabularnewline
			\hline 
			\hline
			1 & 120   & 1864   & 1.29986350 & 257.565054 &       -  \tabularnewline \hline
			2 & 960   & 7356   & 1.30100271 & 238.189022 & 19.37603 \tabularnewline \hline
			3 & 3240  & 18768  & 1.30099274 & 223.291696 & 14.89733 \tabularnewline \hline
			4 & 7680  & 38260  & 1.30100097 & 218.937003 &  4.35469 \tabularnewline \hline
			5 & 15000 & 67992  & 1.30100587 & 217.083298 &  1.85370 \tabularnewline \hline
			6 & 25920 & 110124 & 1.30100827 & 216.105059 &  0.97824 \tabularnewline \hline
		\end{tabular}
  \caption{{Detuning values for the 1-cell TESLA cavity.}\label{tab:FrequencyShift}}
\end{center}
\end{table}

\begin{figure}
\begin{center}
\includegraphics[width=.9\textwidth]{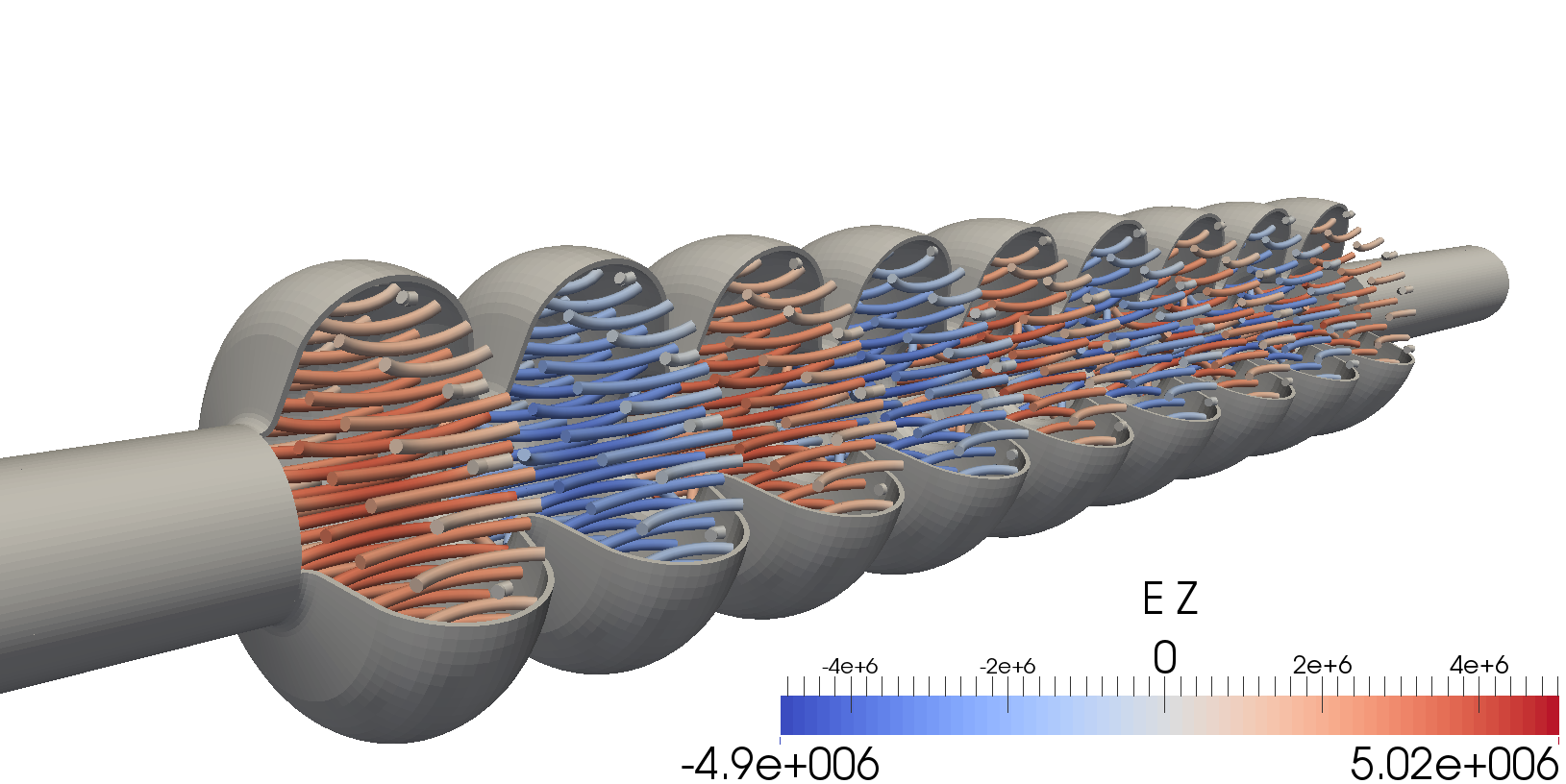}
\caption{The fundamental $TM_{010},\pi$ mode for the 9-cell TESLA cavity.\label{fig:Ez9cell}}
\end{center}
\end{figure}

Starting from \cite{Tesla2000}, the geometry for full 9-cell TESLA cavity has been created. With respect to the single cavity, one has to take into account that, due to the coupling between the different cells, the fundamental mode splits itself into 9 different modes with similar frequencies giving rise to the so-called fundamental passband. The results for these eigenfrequencies are shown in Table~\ref{tab:Frequency9cell}: the accelerating mode is the $\pi$ mode at 1.3~GHz. The $z$ component for the electrical field of the $\pi$ modes is depicted in Fig.~\ref{fig:Ez9cell}.

\begin{table}
\begin{center}
		\begin{tabular}{|c|c|}
			\hline
			Mode & Frequency [GHz]\tabularnewline
			\hline 
			\hline
			1 & 1.276335705889215\tabularnewline
			\hline 
			2 & 1.278421359483793\tabularnewline
			\hline 
			3 & 1.281632725459760\tabularnewline
			\hline
			4 & 1.285597822640840\tabularnewline
			\hline
			5 & 1.289849369271624\tabularnewline
			\hline
			6 & 1.293875642584478\tabularnewline
			\hline
			7 & 1.297181927064266\tabularnewline
			\hline
			8 & 1.299363801453597\tabularnewline
			\hline
			9 & 1.300002415591750\tabularnewline
			\hline
		\end{tabular}
\caption{Frequencies for the first passband of the 9-cell TESLA cavity (151888 DoFs). \label{tab:Frequency9cell}}
\end{center}
\end{table}

\section{Conclusions}
Low order Finite Element Methods may fail to achieve a sufficient accuracy for calculating Lorentz detuning in superconducting accelerator cavities. {This could be alleviated using software where curved elements and methods such as sensitivity analysis can be exploited. Alternatively, this work proposes Isogeometric Analysis as a solution for Lorentz detuning simulation since it naturally comprehend a better representation of the curved cavity walls and a natural way for treating mechanical deformations within the electromagnetic eigenvalue problem, without loss of geometric accuracy. The results show that the Isogeometric method succeeds in obtaining reliable results for the frequency shifts. Furthermore, the tests performed indicate a higher computational efficiency regardless of the different properties of the matrices.} 

\section*{Acknowledgements}

	The authors would like to thank CST AG, i.e. Stefan Reitzinger and Rodrigo Enjiu, for providing support and the FE matrices for the comparison.	
	
	This work is supported by the 'Excellence Initiative' of the German Federal and State Governments and the Graduate School of Computational Engineering at Technische Universit\"at Darmstadt.
	
	Carlo de Falco's work is partially funded by the `Start-up Packages and PhD Program project', co-funded by Regione Lombardia through the `Fondo per lo sviluppo e la coesione 2007-2013', formerly FAS program.

\end{document}